\documentclass{siamart0516}

\usepackage{microtype}

\usepackage{newunicodechar}
\newunicodechar{α}{\ensuremath \alpha}
\newunicodechar{β}{\ensuremath \beta}

\usepackage{amsmath}
\DeclareMathOperator{\E}{E}
\DeclareMathOperator{\Var}{Var}

\usepackage{listings}
\lstdefinelanguage{julia}{
  basicstyle=\small\ttfamily,
  showspaces=false,
  showstringspaces=false,
  keywordstyle={\textbf},
  morekeywords={if,else,elseif,while,for,begin,end,quote,try,catch,return,local,abstract,function,generated,macro,ccall,finally,typealias,break,continue,type,global,module,using,import,export,const,let,bitstype,do,in,baremodule,importall,immutable},
  escapeinside={~}{~},
  morecomment=[l]{\#},
  commentstyle={},
  morestring=[b]",
}
\usepackage{algorithm}
\usepackage{algorithmic}
\usepackage{pgfplots}
\usepackage{subcaption}

\usepackage{pifont}
\newcommand{\cmark}{\ding{51}}%
\newcommand{\xmark}{\ding{55}}%

\title{Fast computation of the principal components of genotype matrices in Julia
    \thanks{This
        work was supported by the Intel Science and Technology Center for Big Data.
	}}

\author{%
    Jiahao Chen
    \thanks{Computer Science and Artificial Intelligence Laboratory,
           Massachusetts Institute of Technology,
           Cambridge, Massachusetts 02139.
           Current Address: Capital One Financial,
           11 West 19th Street,
           New York, New York 10011.
           ({\tt jiahao.chen@capitalone.com})}
    \and
    Andreas Noack
    \thanks{Computer Science and Artificial Intelligence Laboratory,
            Massachusetts Institute of Technology,
            Cambridge, Massachusetts 02139.
            Current address: Julia Computing,
            114 Western Ave, Allston, Massachusetts 02134.
            ({\tt andreasnoackjensen@gmail.com})}
    \and
    Alan Edelman
    \thanks{Department of Mathematics and Computer Science and Artificial Intelligence Laboratory,
            Massachusetts Institute of Technology,
            Cambridge, Massachusetts 02139 ({\tt edelman@mit.edu})}
}

\begin{document}

\maketitle

\tableofcontents

\begin{abstract}
Finding the largest few principal components of a matrix of genetic data
is a common task in genome-wide association studies (GWASs), both for
dimensionality reduction and for identifying unwanted factors of
variation.  We describe a simple random matrix model for matrices that
arise in GWASs, showing that the singular values have a bulk behavior that
obeys a Marchenko-Pastur distributed with a handful of large outliers.
We also implement Golub-Kahan-Lanczos (GKL) bidiagonalization in the Julia
programming language, providing thick restarting and a choice between full
and partial reorthogonalization strategies to control numerical roundoff.
Our implementation of GKL bidiagonalization is up to 36 times faster than
software tools used commonly in genomics data analysis for computing principal
components, such as EIGENSOFT and FlashPCA, which use dense LAPACK routines
and randomized subspace iteration respectively.
\end{abstract}

\begin{keywords}
    singular value decomposition,
    principal components analysis,
    genome-wide association studies,
    statistical genetics,
    Lanczos bidiagonalization,
    Julia programming language,
    subspace iteration
\end{keywords}

\begin{AMS}
    65F15, 97N80
\end{AMS}

\pagestyle{myheadings}
\thispagestyle{plain}
\markboth{J. CHEN, A. NOACK AND A. EDELMAN}{Principal components of genotype matrices in Julia}

\section{Principal components of genomics data}

Personalized medicine or precision medicine is a growing movement to tailor treatments of disease
to an individual's sensitivities to treatment, allergies, or other genetic
predispositions, using all available data about an individual~\cite{Desmond2012}.
Developers of personalized medical treatments are therefore interested in how an
individual's genome, both in isolation and within the context of the wider human
population, can be used to predict desired clinical outcomes (known as
comorbidities or phenotypes)~\cite[Ch. 8]{Laird2011}.
Genome-wide association studies (GWASs) are a new and popular technique for
studying the significance of human genome data, by studying the associate genotype
variation with phenotype variations in outcome variables e.g.\ the clinical
observation of a disease.

The genome data used in
GWASs are are often encoded in a matrix expressing the number of mutations
from a reference genome, which we will refer to as a \textit{genotype matrix}.
By convention, the genotype matrix is indexed by patients (or other test
subjects) on the rows and gene markers on the columns which represent some
coordinate or locus within the human genome. A common example of gene markers are
single nucleotide polymorphisms (SNPs), which represent gene positions with
pointwise mutations of interest that express variation in the genotype across
humans. Oftentimes, the explanatory data in a GWAS are simply called SNP data.
Since most human cells are at most diploid, i.e.\ have two sets of chromosomes,
the matrix elements can only be 0, 1, or 2 (or missing).

There are two major confounding sources of variation which are considered in
the analysis of human genome data, each of which have significance for the
spectral properties of the matrix:

\paragraph{Population stratification/admixing}
Population stratification is the phenomenon of common genetic variation within
mutually exclusive subpopulations defined along racial, ethnic or geographic
divisions~\cite{Pritchard1999,Cardon2003}. (Admixture models relax the mutual
exclusion constraint~\cite{Devlin1999,Sankararaman2008}.)
In linear algebraic terms, the genotype matrix will have a low rank component
with large singular values.

\paragraph{Cryptic relatedness}
Sometimes called kinship or inbreeding, cryptic relatedness is an increase in
sampling bias in the columns (human genomes) produced by having common ancestors,
thus increasing the nominally observed frequency of certain
mutations~\cite{Voight2005,Astle2009}. Relatedness is usually detected and removed in
a separate preprocessing step, but it is not always possible to remove
fully~\cite{PLINK}. Any remaining related samples will
result in (near) linear dependencies in the rows of the genotype matrix,
leading to the presence of several singular values that are very small or zero.

Principal component analysis (PCA) was historically first used as a dimensionality
reduction technique to summarize the variation in the human genes and study its
implications for human evolution in relation to other factors such as
geography and history~\cite{Menozzi1978,Cavalli1994,Novembre2008}.
However, we will focus on the more modern use of principal components (PCs) to
represent the confounding effects of population substructure in the statistical
modeling of GWASs~\cite{Chen2003,Patterson2006,Price2006,Zhu2002,Zhang2003}.

Genomics matrices form an interesting use case for the classical techniques of
numerical linear algebra, as the amount of sequenced genome data grows
exponentially~\cite{Stephens2015}. As the price of sequencing genome data declines
rapidly, genomics studies involving hundreds of thousands of individuals (columns)
are already commonplace today, with order of magnitude growth expected within the next year or two. Therefore, genomics researchers will require access to the best
available algorithms for parallel computing and numerical linear algebra to
handle the increasing demands of data processing and dimensionality reduction.

\subsection{The statistical significance of principal components}

The main statistical tool used in GWAS is regression, using some model that
associates genotype variation with phenotype variation.
While correlation does not imply causation in and of itself, the central dogma
of molecular biology states that causality flows from genetic data in DNA and RNA
to phenotype data in expressed proteins~\cite{Crick1970}. Consequently, correlations between genotypes and phenotypes could in theory have causal significance.
The \emph{linear} regression model is one the simplest useful models, and can be
motivated in several different ways. One way is in terms of least squares
minimization to find the coefficients that minimize the sum of squared distances
between a hyperplane and the observations. Another way to formulate the problem
as finding a projection of the vector of observations down onto the space
spanned by the vectors of explanatory variables.
One caveat in statistical studies is the assumption that the observations are
randomly sampled with replacement, and hence that the observations can be
assumed to be independent. In the context of statistical genomics, the
independence assumption is one of several bundled into the principle of
Hardy-Weinberg equilibrium~\cite{Hardy1908,Weinberg1908}, which is commonly
assumed in statistical genetics~\cite{Laird2011}. However, caution must be taken
to remove possible sampling bias due to the collection of genetic data from
patients from a single hospital or hospital network, on top of sampling bias
introduced by not treating the presence of population substructure.

\subsection{The linear regression model}

The statistical theory of regressing phenotypes against genotypes is best expressed
in terms from conditional expectations. If we for individual $i$ denote the genotype measurement by $x_i$ and the phenotype measurement by $y_i$, the conditional
expectation of interest can be written as $\E(y_i|x_i)=\beta_0 + \beta_1 x_i$.
A popular formulation of the this model is
\begin{equation*}
    y_i = \beta_0 + x_i \beta_1 + \varepsilon_i\quad i=1,\dots,n
\end{equation*}
where $n$ is the number of observations which in this case would be the number of individuals for which we have genotype data.
The variable $\varepsilon_i$ is called the error term and must satisfy $\E(\varepsilon_i|x_i)=0$.
More generally, $\varepsilon_i$ is the conditional distribution $y_i - \beta_0 - x_i\beta_1|x_i$.

The multivariate expression can be written conveniently in matrix form

\begin{equation}
\label{eq:linreg}
    y = X\beta + \varepsilon
\end{equation}
with
\begin{equation*}
    X = \begin{pmatrix}
        1 & x_1\\
        \vdots & \vdots \\
        1 & x_n
    \end{pmatrix}
\end{equation*}
and in this notation, the well-known least squares estimator for the coefficient
$\beta$ can be written as
\begin{equation}
\label{eq:ols}
    \hat{\beta} = (X^TX)^{-1}X^Ty.
\end{equation}

The conditional probability treatment demonstrates that when $(x_i,y_i)$ pairs
are considered to be random variables, then $\hat{\beta}$ is also a random
variable, with its mean and variance quantifying uncertainty about the least
squares solution.
First, notice that the least squares estimator \eqref{eq:ols} can be written
\begin{equation*}
    \hat{\beta} = \beta + (X^TX)^{-1}X^T\varepsilon
\end{equation*}
and the expected value of the estimator is
\begin{equation}
    \E(\hat{\beta}|X) = \E(\beta + (X^TX)^{-1}X^T\varepsilon|X) =
        \beta + (X^TX)^{-1}X^T\E(\varepsilon|X) = \beta
\end{equation}
which means that the estimator is \emph{unbiased}. Statisticians are interested in the variability of $\hat{\beta}$ under changes to the data that could be considered small errors. The most common measurement for the variability of an estimator is the (conditional) variance, i.e.
\begin{equation*}
    \Var(\hat{\beta}|X) = \Var(\beta + (X^TX)^{-1}X^T\varepsilon|X) =
        (X^TX)^{-1}X^T\Var(\varepsilon|X)X(X^TX)^{-1}.
\end{equation*}

This shows that variance of $\hat{\beta}$ depends on the (conditional) variance of $\varepsilon$, which has not been discussed yet. In classical treatments of the linear regression model it is typically assumed that, conditionally on $x_i$, the $y_i$s are independent and the have the same variance which is the same as $\Var(\varepsilon|X)=\sigma^2 I$ for some unknown scalar $\sigma^2$. Under this assumption, the variance of $\hat{\beta}$ reduces to $\sigma^2 (X^TX)^{-1}$. The magnitude of this quantity is unknown because $\sigma^2$ is an unknown parameter but $\sigma^2$ can be \emph{estimated} from the data. The usual estimator is
$\hat{\sigma}^2=\frac{1}{n} \|\hat{\varepsilon}\|^2$
where $\hat{\varepsilon} = y - X\hat{\beta}$.
This leads to the estimate of the (conditional) variance of the estimator $\widehat{\Var(\hat{\beta}|X)} = \hat{\sigma}^2(X^TX)^{-1}$.

The independence assumption is often used in statistics and can be justified from an assumption of random sampling with replacement. In studies where data is passively collected, this might not be a reasonable assumption as explained in the previous section. Non-random sampling might lead to correlation between the phenotypes even after conditioning on the genotypes. In consequence of that, $\Var(\varepsilon|X)$ will no longer be diagonal, but have some general positive definite structure $\Sigma$ and $\Var(\hat{\beta}|X)=(X^TX)^{-1}X^T\Sigma X(X^TX)^{-1}$. Since $\Sigma$ in general consists of $\frac{n(n+1)}{2}$ parameters, it cannot be estimated consistently.

In order to analyze the problem with correlated observations, it is convenient to decompose the error into a part that contains the cross-individual correlation and a part that is diagonal and therefore only describes the variance for each individual. This may be written as
\begin{equation*}
    y_i = \beta_0 + x_i\beta_1 + \eta_i + \xi_i
\end{equation*}
where is assumed that $\Var(\eta|X)=\Sigma_\eta$ and $\Var(\xi|X)=\sigma^2_\xi I$. Furthermore, it is assumes that the two error terms are independent.

\label{sec:regress-correction}
\subsubsection{Fixed effect estimation}

One way to produce a reliable estimate of the variance of $\hat{\beta}$ is to come up with a set of variables $z_1,\dots,z_k$ that proxies the correlation between the observations, i.e.\ $\eta=Z\gamma$. By simply including the variables $z_1,\dots,z_k$ in the regression model, it possible to remove the correlation which distorts the variance estimate for $\hat{\beta}$. For the regression $y|X,Z$, we get the least squares estimator
\begin{equation*}
    \begin{pmatrix}\hat{\beta} \\ \hat{\gamma}\end{pmatrix} =
    \begin{pmatrix}
        X^TX & X^TZ \\ Z^TX & Z^TZ
    \end{pmatrix}^{-1}\!\!
    \begin{pmatrix}
        X^T \\ Z^T
    \end{pmatrix}
    (X\beta + Z\gamma + \xi)=
    \theta +
            \begin{pmatrix}
        X^TX & X^TZ \\ Z^TX & Z^TZ
    \end{pmatrix}^{-1}\!\!
    \begin{pmatrix}
        X^T \\ Z^T
    \end{pmatrix}\xi,
\end{equation*}
which has variance
\begin{equation*}
    \Var\left(\begin{pmatrix}\hat{\beta} \\ \hat{\gamma}\end{pmatrix}|X,Z\right)=\sigma_\xi^2
        \begin{pmatrix} X^TX & X^TZ \\ Z^TX & Z^TZ \end{pmatrix}^{-1}.
\end{equation*}

In many applications, a few principal components of the covariance matrix of the complete SNP data set is used as a proxy for the correlation between individuals. Computing principal components is therefore often a first step in analyzing GWASs.

\subsection{Software for computing the principal components of genomics data}

The software stack for GWAS is generally based on command line tools written completely in C/C++. Not only is the core computational algorithm written in C/C++ but also much of the pre- and postprocessing of the data. The data sets can be large but and the computations at times heavy but it has been a surprise to learn the extend to which analyses are carried out directly from the command line instead of using higher level languages like MATLAB, R, or Python. This choice seems to limit the tools available to the analysts because, unless he is a C/C++ programmer, the programmer is restricted to the set of options included in the command line tool.

Two major packages exist for computing the PCA in the GWAS software stack. The package \emph{EIGENSOFT} accompanied~\cite{Patterson2006} which popularized the use of PCA in GWAS. In the original version of the package, the routine \texttt{smartpca} for computing the PCA of a SNP matrix was based on an eigenvalue solver from LAPACK. In consequence, all the eigenvalues and vectors of the SNP matrix were computed even though only a few of them were used as principal components. Computing the full decomposition is inefficient and as the number of available samples has grown over the year, this approach has become impractical.

More recently, \emph{FlashPCA}~\cite{abraham2014fast} has emerged as a potentially faster alternative to EIGENSOFT's \texttt{smartpca}. The PCA routine is based on a truncated SVD algorithm described in~\cite{halko2011finding}. More specifically, FlashPCA uses a subspace iteration scheme with either column-wise scaling (FlashPCA 1) or orthogonalization by QR (FlashPCA 2) in each iteration. The convergence criterion is based on the average of squared element-wise distance between the the bases for iteration $i-1$ and $i$. The QR orthogonalization step in the implementation of FlashPCA routine deviates from the algorithm in described~\cite{abraham2014fast} which only normalizes column-wise. Our conjecture was that this change was made to avoid loss of orthogonality in the subspace basis and the author of the package has confirmed this. In consequence, the timings in~\cite{abraham2014fast} do not correspond to the performance of the software run with default settings since the QR orthogonalization is much slower than the column-wise normalization. Furthermore, a degenerate basis might also converge much faster because it eventually converges to the single largest eigenvalue.

Table~\ref{tab:sw} lists some software packages providing eigenvalue or singular
value computations useful for PCA. It may surprise readers that well-established
packages in numerical linear algebra are rarely used in genomics, considering
that libraries like ARPACK and PROPACK have convenient wrapper libraries in both
R, Python, and \textsc{MATLAB}. This phenomenon might be explained by the pronounced use
of C/C++ in statistical genetics, where calling ARPACK and PROPACK are relatively
more demanding, combined with the fact that iterative methods are traditionally
not part of the curriculum in statistical genomics.

\begin{table}
    \centering
    \caption{List of software used for principal components analysis in statistical
    genomics. Notably, common packages known to numerical linear algebraists are
    rarely, if ever, used.
    $^*$: available in EIGENSOFT~\cite{Patterson2006}.
    $^\dagger$: uses dense LAPACK routines.
    \label{tab:sw}}
    \begin{tabular}{cccc}
       \hline
        Software         & Reference & Algorithm                              & Used in genomics \\ \hline
        smartpca$^{*\dagger}$ & \cite{Patterson2006} & Householder bidiagonalization & \cmark \\ FastPCA$^*$ & \cite{Galinsky2016} & Subspace iteration        & \cmark \\
        FlashPCA         & \cite{abraham2014fast} & Subspace iteration      & \cmark \\ \hline
        ARPACK           & \cite{Lehoucq1996} & Lanczos tridiagonalization  & \xmark \\
        PROPACK          & \cite{Larsen1998} & Lanczos bidiagonalization    & \xmark \\
        SLEPc            & \cite{Hernandez2008} & Lanczos bidiagonalization & \xmark \\
        Anasazi          & \cite{Baker2009} & Lanczos bidiagonalization     & \xmark \\ \hline
    \end{tabular}
\end{table}

\label{sec:model}
\section{A simple model for genomics data}

In this section we present a very simple random model which accurately mimics
the spectral features observed in real data matrices.
We hypothesize that the spectral properties of human patient genotype data can
be modeled the Julia code in Algorithm~\ref{code:model}, which captures the
confounding effects of population admixture and cryptic relatedness.
The former is modeled by setting randomly select subblocks to the same value,
whereas the latter is modeled by duplicating rows, thus purposely introducing
linear dependence into the row space. This model, while simplistic, can be
tuned to reproduce the scree plot observed in empirical data matrices and we
expect that such models may be of interest to researchers developing numerical
algorithms that lack access to actual data, which are often access-restricted
due to clinical privacy.

\begin{algorithm}
\caption{A simple model for human genotype data matrices in Julia
\label{code:model}}
\begin{lstlisting}
"""
Inputs:
- m: number of rows (gene markers)
- n: number of columns (patients)
- r: number of subblocks to model population admixing
- nsignal: number of entries to represent signal
- rkins: fraction of columns to duplicate

Output:
- A: a dense matrix of size m x n with matrix elements 0, 1 or 2
"""
function model(m, n, r, nsignal, rkins)
    A = zeros(m, n)
    #Model population admixing
    #by randomly setting a subblock to the same value, k
    for i=1:r
        k = rand(0:2)
        r1 = randrange(m)
        r2 = randrange(n)
        A[r1, r2] = k
    end

    #Model signal
    for i=1:nsignal
        A[rand(1:m), rand(1:n)] = rand(0:2)
    end

    #Model kinship by duplicating rows
    nkins = round(Int, rkins*m)
    for i=1:nkins
        A[rand(1:m), :] = A[rand(1:m), :]
    end
    return A
end

function randrange(n)
    i1 = rand(1:n)
    i2 = rand(1:n)
    if i1 > i2
        return i2:i1
    else
        return i1:i2
    end
end
\end{lstlisting}
\end{algorithm}

Algorithm~\ref{code:model} describes a \verb|model| function, which creates a
dense, synthetic data matrix of size $m\times n$. The synthetic data is
generated in three steps. First, randomly select $r$ rectangular submatrices
(which may overlap) and set the elements of each submatrix to the same value.
This process simulates crudely the effects of population admixture, where
each subpopulation has a common block of mutations that vary together, and the
possibility of overlap resembles the effect of mixing different subpopulations
together. (This step uses the auxiliary \verb|randrange| function, which returns
a valid Julia range that is a subinterval of the range \verb|1:n|.)
Second, model the genotype of interest by randomly setting $nsignal$ matrix
elements randomly to 0, 1 or 2.
Third, simulate the effects of kinship by choosing a fraction $rkins$ of the
rows to duplicate.

\begin{figure}

\includegraphics[width=0.45\textwidth]{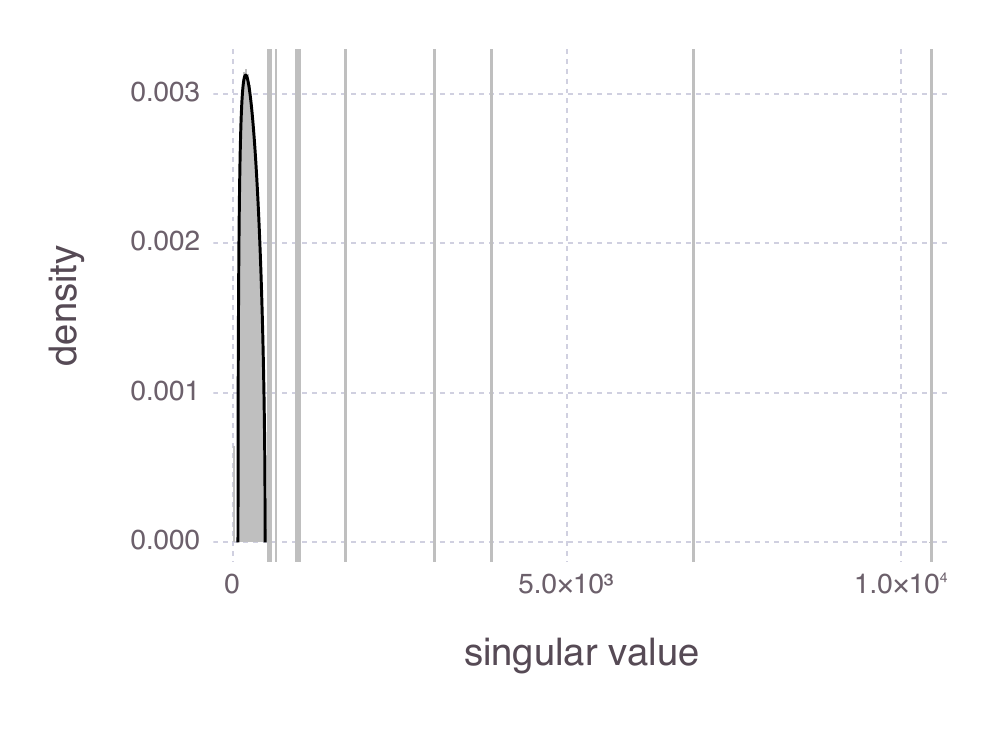}
\includegraphics[width=0.45\textwidth]{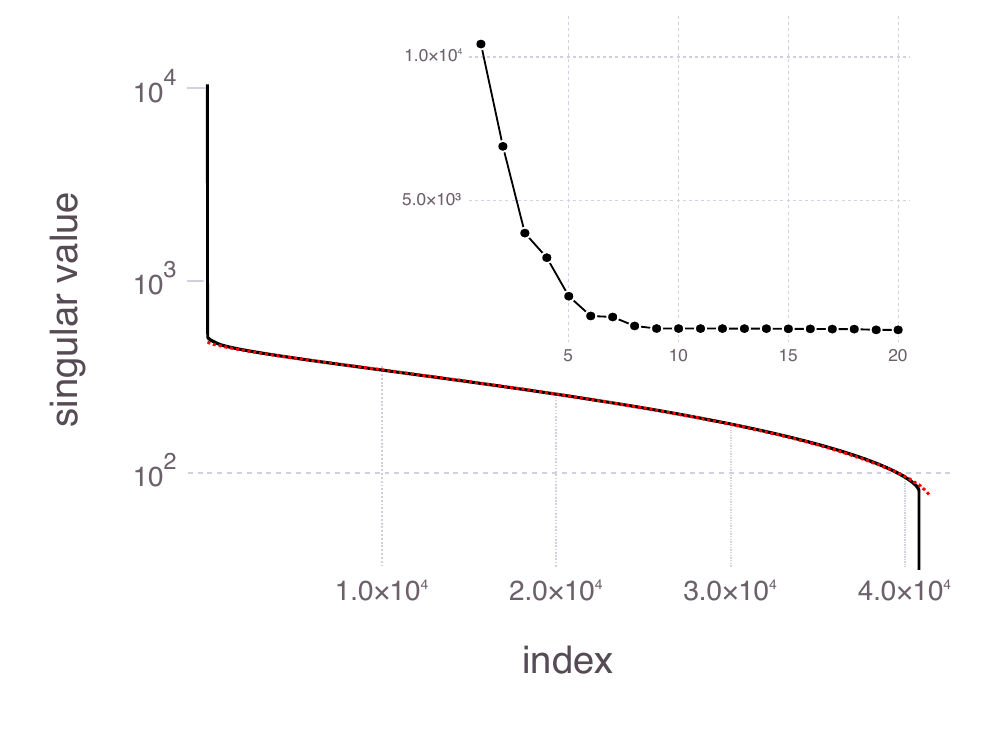}

\caption{Left: Histogram of the singular values of a synthetic genomics data
matrix of size $41505\times81700$ (grey bars) generated with the code in
Algorithm~\ref{code:model}, overlaid with the
Marchenko-Pastur law \eqref{eq:mplaw-sv} for $\rho = 1.86$ (black line).
Light gray vertical lines show the
presence of 18 outlying singular values whose magnitudes exceed
$1.1\sigma_+ \approx 529.7$.
Right: Scree plot of the singular values of the same matrix
(black solid lines), showing the presence of a large, low rank portion
of approximately rank 10 (inset), and an asymptotic convergence to the same
Marchenko-Pastur law (red dotted line).
\label{fig:empirical-spectrum}
}
\end{figure}

The model is clearly a very crude approximation to genotype data with the
confounding effects of population substructure. There is little overt control
over the admixture process, the precise distribution over matrix elements, and
all related patients are assumed to have identical genotypes. Nevertheless, the
model generates realistic distributions of singular values which mimic closely
what we have observed in that of real world genotype matrices.
Figure~\ref{fig:empirical-spectrum} shows the distribution of singular values
generated by our model with parameters $m=41505$, $n=81700$, $r=10$,
$nsignal=mn$, $rkins=0.017$.
There are a handful of ($\approx r$) large singular values,
while the rest follow a bulk distribution from random matrix theory known as the
Marchenko-Pastur law~\cite{Marchenko1967} with parameter $\rho = 1.86$.

The Marchenko-Pastur law describes the distribution of governs the eigenvalues
of a random covariance matrix $Y=XX^T$ formed from a data matrix $X$
of iid elements with mean 0 and finite variance $\sigma^2$.
Let $\rho$ be the ratio of the number of rows of $X$ to the number of columns of
$X$. Then, the nonzero eigenvalues follow the distribution
\begin{equation}
    \label{eq:mplaw-ev}
    p_e(\xi) = \frac {\sqrt{(\lambda_+-\xi)(\xi-\lambda_-)}} {2 \pi \sigma^2 \lambda x}
\end{equation}
where
$\lambda_+ = \sigma^2(1+\sqrt{\rho})^2$ and
$\lambda_- = \sigma^2(1-\sqrt{\rho})^2$.

When written in terms of the probability density of the singular values of $X$,
the law reads
\begin{equation}
    \label{eq:mplaw-sv}
    p_s(x) = \frac {\sqrt{(\sigma_+^2-x^2)(x^2-\sigma_-^2)}} {\pi \sigma^2 \min(1, \lambda) x}
\end{equation}
where
$\sigma_+ = \sqrt{\lambda_+} = \sigma(1+\sqrt{\rho})$ and
$\sigma_- = \sqrt{\lambda_-} = \sigma\vert1-\sqrt{\rho}\vert$.

\begin{figure}
\caption{Marchenko-Pastur law for $\rho=1.5$ (black lines) for the densities of
nonzero eigenvalues of a random covariance matrix ($XX^{T}$) and the singular
values of $X$ for iid matrix elements with mean 0 and variance 1
(black lines), corresponding to \eqref{eq:mplaw-ev} and \eqref{eq:mplaw-sv}
respectively.
Shown for comparison are corresponding histograms (grey
bars) of numerically sampled eigenvalues and singular vectors from
a numerically sampled random matrix $X$ of size $1000\times667$,
with iid Gaussian entries of mean 0 and variance 1.
\label{fig:mplaw}
}

\includegraphics[width=0.45\textwidth]{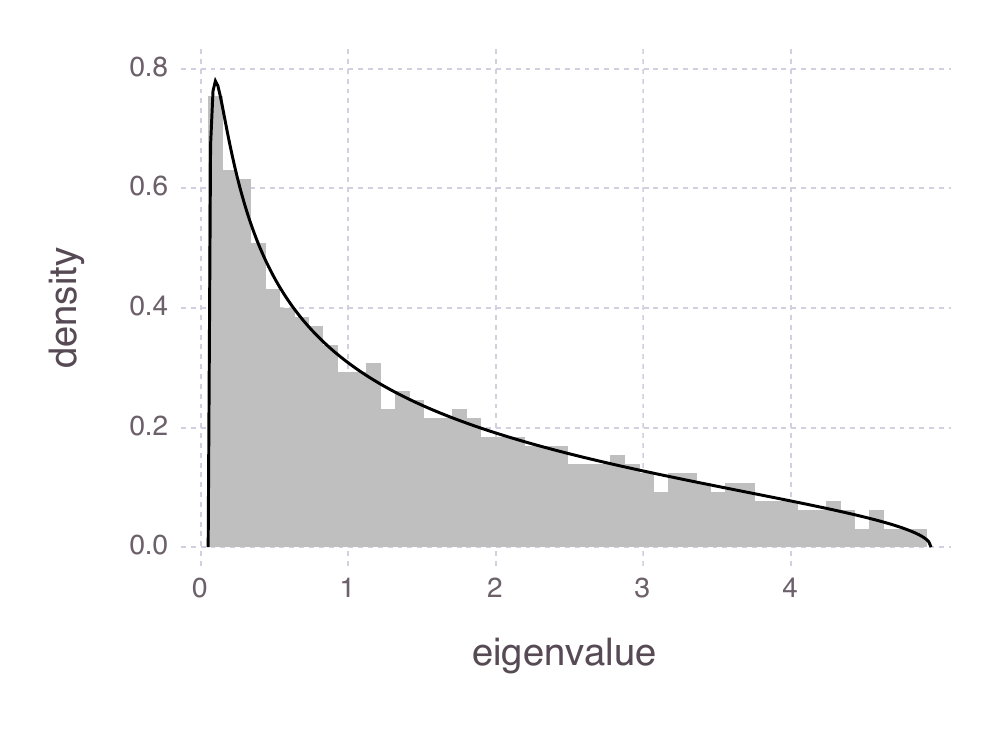}
\includegraphics[width=0.45\textwidth]{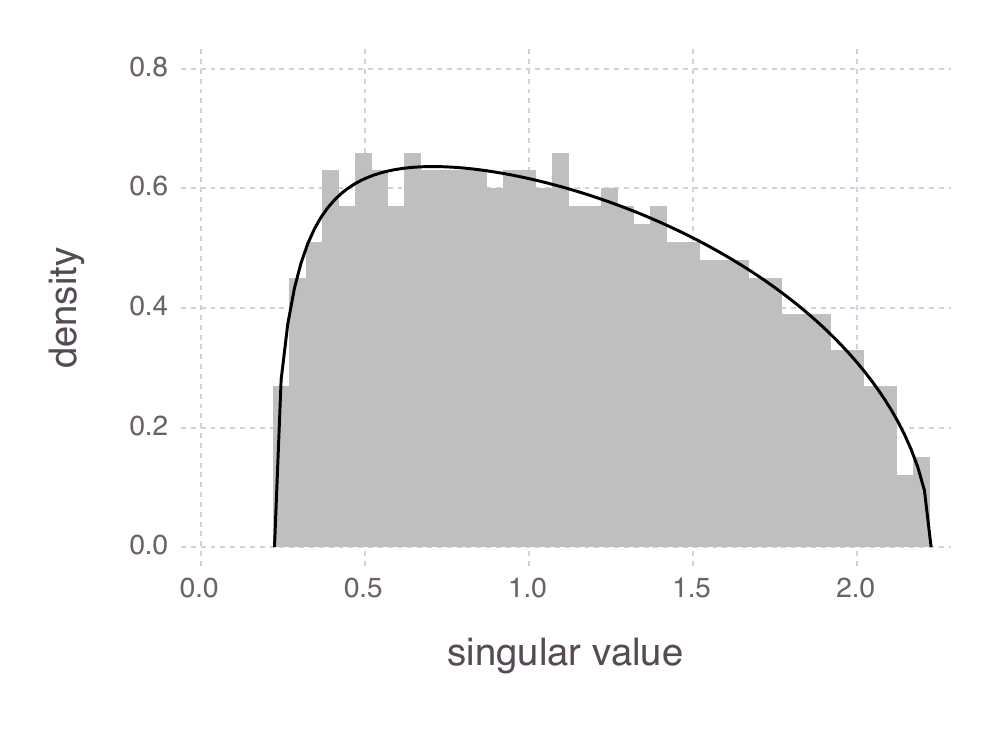}
\end{figure}

Figure~\ref{fig:mplaw} shows typical density plots for random eigenvalues
and singular values for $\rho=1.5$.

It is worth noting that random matrix theory had been previously introduced in
the theoretical analysis of principal components of genotype matrices.
\cite{Patterson2006} proposed a hypothesis test that computed principal components
should correspond to eigenvalues that were different from those expected from a
pure random covariance matrix,
and comparing in particular the largest eigenvalue against one randomly sampled
from the Tracy-Widom distribution~\cite{Tracy1993,Tracy1994}.
This analysis only shows that the matrix is not iid. In contrast, we show here
an explicit construction of a random matrix, whose elements are not iid, that
can generate a realistic spectrum of singular values that consists of several
large outliers and a bulk distribution that empirically satisfies the
Marchenko-Pastur law, albeit with a modified parameter $\rho = 1.86$ as opposed
to the value $81700/41505=1.97$ which would be expected from taking the ratio of
the number of rows to the number of columns.

\section{Algorithms for PCA}

The discussion in Section~\ref{sec:model} demonstrates that the
confounding effects of population substructure can and does produce a low
rank structure in the top singular vectors (i.e.\ singular triples
corresponding to the largest singular values), which can be captured even
in the very crude random matrix model of Algorithm~\ref{code:model}. The
top few principal components,
which by construction capture the largest components of the variability, are
good candidates for modeling the unwanted variation as described in
Section~\ref{sec:regress-correction} and have been used in the statistical
genetics community for this
purpose~\cite{Chen2003,Patterson2006,Price2006,Zhu2002,Zhang2003}.

Iterative eigenvalue (or singular value)
methods~\cite{Bai2000} are therefore computationally efficient choices for
determining these principal components, as only a handful of them are needed.
However, none of the classical methods known in numerical linear algebra
(apart from subspace iteration) is implemented in commonly used software packages
for PCAs in genomics, as shown in Table~\ref{tab:sw}.
Notably absent is any Lanczos-based bidiagonalization method.
We have therefore implemented the Golub-Kahan-Lanczos bidiagonalization~\cite{Golub1965}
method in pure Julia\cite{Bezanson2012,Bezanson2015}.
We have incorporated several of the best available numerical features, such as:

\paragraph{Thick restarting}
To control the memory usage and accumulation of
roundoff error, we also include an implementation of the thick restart
strategy~\cite{Wu2000,Stewart2001}, offering restarts using either
ordinary Ritz values~\cite{Wu2000} or harmonic Ritz
values~\cite{Baglama2005}. The thick restart variant is becoming
increasingly popular, being available not only in
SLEPc~\cite{Hernandez2008}, but also in R as the IRLBA
package~\cite{cran-irlba}, and has also become the new algorithm for
\verb|svds| in \textsc{MATLAB} R2016a (which also offers partial
reorthogonalization)~\cite{matlab-r2016a}. For the purposes of computing
the top principal components, it suffices to use thick restarts using
ordinary Ritz values.

\paragraph{Choice of reorthogonalization strategy}
We offer users the choice of partial
reorthogonalization~\cite{Simon1984,Larsen1998} or full
reorthogonalization using doubly reorthogonalized classical Gram-Schmidt.
By default, the implementation uses an adaptive threshold for determining
when the second reorthogonalization is necessary, based on the expected
number of digits lost to catastrophic
cancellation~\cite{Daniel1976,Bjorck2015}.
We also use the adaptive one-sided reorthogonalization strategy on either
the left or right singular vectors (whichever is smaller), unless our
estimate of the matrix norm is sufficiently large so that two-sided
reorthogonalization is necessary~\cite{Simon2000}.

\paragraph{Convergence criteria}
We use several different tests for determining when a Ritz value has
converged. At the beginning of the calculation when no other information
is known about the Ritz values, we use the crude estimate on the absolute
error bound on the singular values based on the residual norm computed
from a candidate Ritz value-vector pair~\cite[Ch. 3, \S 53, p.
70]{Wilkinson1965}. However, when the Ritz values become sufficiently
well-separated, more refined estimates can be derived from the
Rayleigh-Ritz properties of the Krylov process~\cite[Ch. 3, \S 54-55, p.
73]{Wilkinson1965}\cite{Yamamoto1980,Ortega1990}.
Experimental facilities are also provided to print and inspect further
convergence information, such as Yamamoto's eigenvector error
bounds~\cite{Yamamoto1980}, Geurt's formula for the componentwise backward
error~\cite{Geurts1982}, Deif's results for \textit{a posteriori} bounds
on eigenpairs and their backward error~\cite{Deif1989}. Interested users
and developers can easily modify the code to implement and inspect yet
other other proposed termination criteria~\cite{Bennani1994}.

\begin{figure}
\caption{Golub-Kahan-Lanczos bidiagonalization in Julia with no restarting and
partial reorthogonalization at a threshold of $\omega = 10^{-8}$, requesting
the top 20 singular values only. Orange vertical lines show when reorthogonalization
was triggered in the computation.
Left: Convergence behavior of the Ritz values.
Right: Computed errors in Ritz values from the residual norm criterion.}
\label{fig:conv-lanczos-tsvd}
\includegraphics[width=0.45\textwidth]{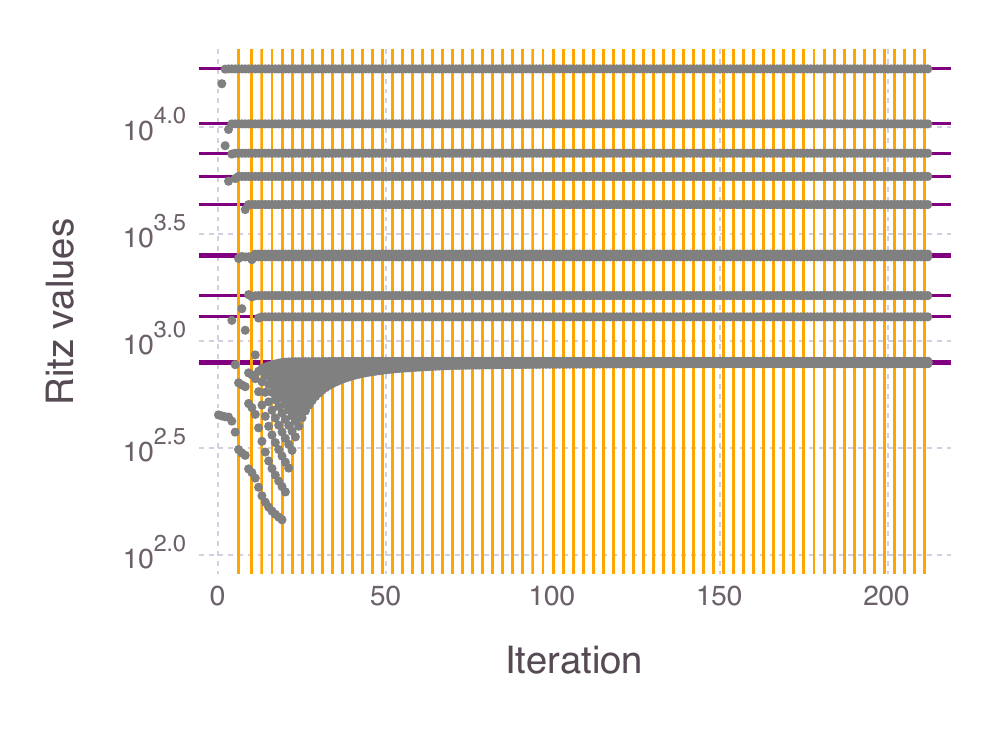}
\includegraphics[width=0.45\textwidth]{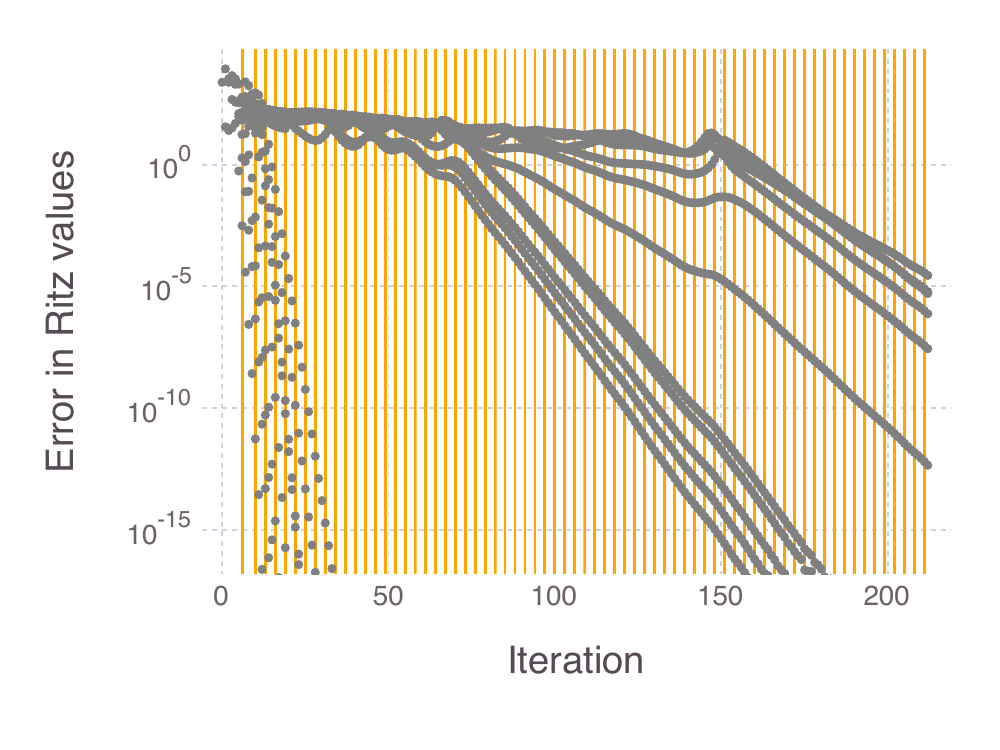}
\end{figure}

\begin{figure}
\caption{Golub-Kahan-Lanczos bidiagonalization in Julia with thick restarting
every 40 microiterations and full reorthogonalization,
requesting the top 20 singular values only.
Gray vertical lines show when restarting was triggered in the computation.
\label{fig:lanczos-tr}
}

\includegraphics[width=0.45\textwidth]{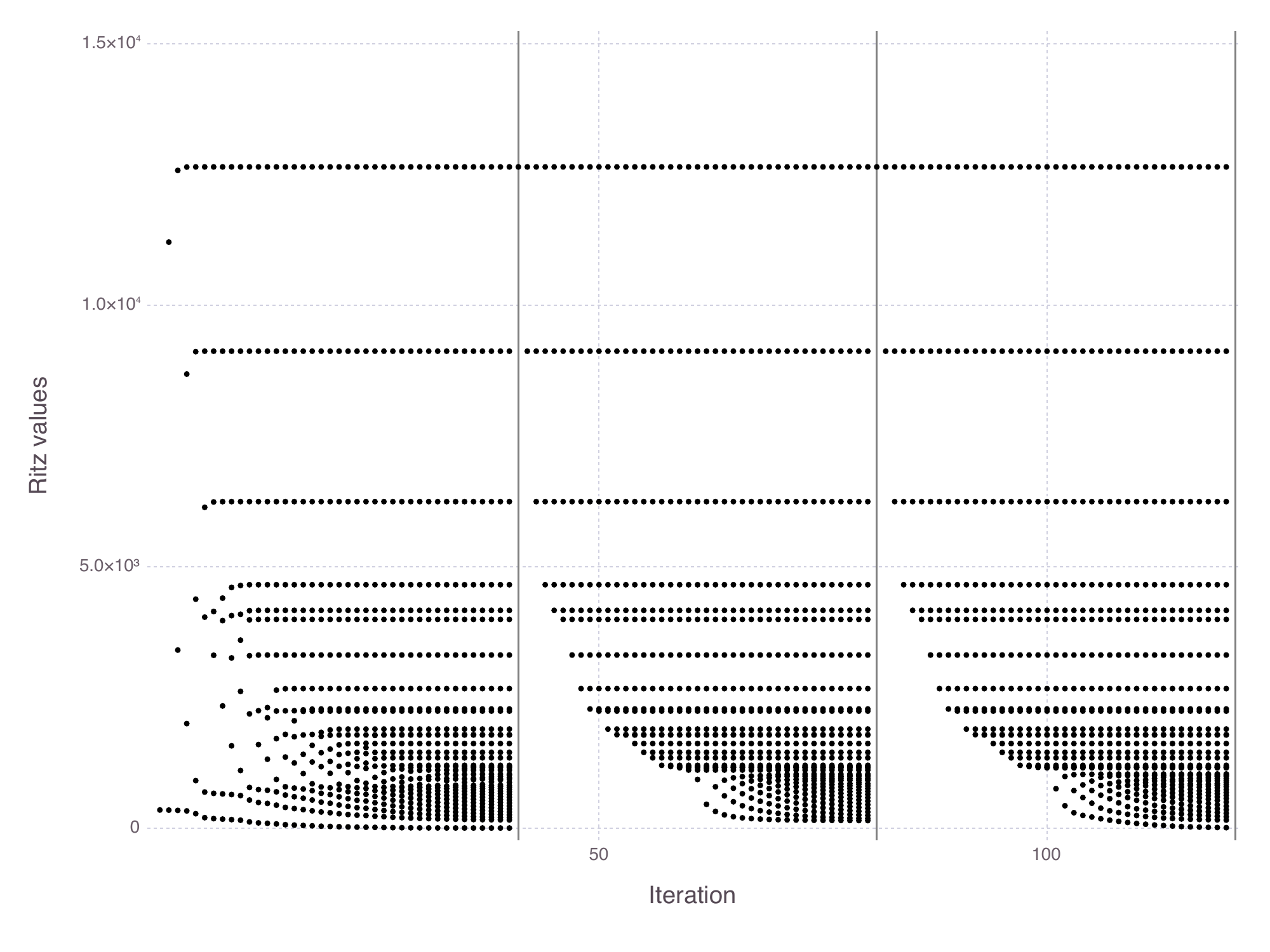}
\end{figure}

The implementation of Lanczos bidiagonalization in Julia allows us to introspect
in great detail into the inner workings of the algorithm.
Figures~\ref{fig:conv-lanczos-tsvd} and \ref{fig:lanczos-tr} show two different
sets of running parameters for our code. The former is run with no restarting
and partial reorthogonalization with threshold $\omega = 10^{-8}$, whereas the
latter uses full reorthogonalization with thick restarting with a maximum
subspace size of 40.

\subsection{Convergence analysis of subspace iteration}

The software in Table~\ref{tab:sw} rely heavily on randomized block power
interactions, which are essentially equivalent to subspace iteration with a
randomized starting subspace.

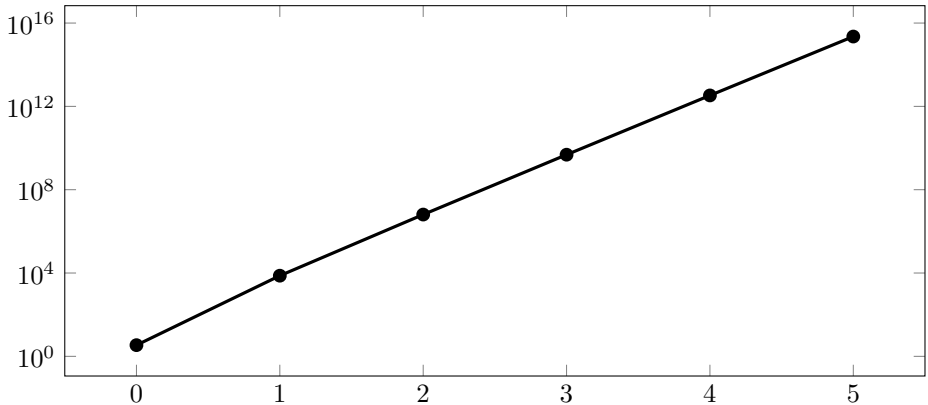
\begin{figure}
    \centering
    \begin{tikzpicture}
        \begin{semilogyaxis}[width=\textwidth, height=0.5\textwidth,
            xtick = {0,1,2,3,4,5}]
            \addplot[very thick, mark=*] file {data/condNumberNorm.txt};
        \end{semilogyaxis}
    \end{tikzpicture}
    \caption{Condition number of the matrix of subspace basis per iteration for
    the published version of FlashPCA~\cite{abraham2014fast}, which normalizes
    each vector but does not reorthogonalize them. Iteration zero is computed
    for the normalized basis after the initial range computation.
    \label{fig:cond}
    }
\end{figure}

FlashPCA also uses an unconventional convergence criterion, namely the matrix
norm of the difference between successive subspace basis vectors,
\begin{equation}
    \label{eq:flashpca-convergence}
    \Delta Y_n = \Vert Y_n - Y_{n-1} \Vert.
\end{equation}
A simple check of the condition number of $Y_n$ with block iteration $n$ shows
that the basis produced by the published version of FlashPCA (which does not
reorthogonalize the basis vectors) rapidly leads to a linearly dependent set of
vectors. As shown in Figure~\ref{fig:cond}, the loss of linear dependence occurs
essentially exponentially, with the inverse condition number reaching machine
epsilon after just five block iterations.
Therefore, we do not recommend the published version of the FlashPCA algorithm
for finding principal components, as it is practically guaranteed to not find
all the principal components requested. We note that the current version of
FlashPCA, which is newer than the published version, reorthogonalizes the
basis vectors by default. In this paper, we refer to the older and newer
versions as FlashPCA 1 and 2 respectively, to avoid ambiguity.

%
%
% \cite[Theorem 4.5.1]{Parlett1998} states, in slightly altered notation, that for
% any scalar $\theta$ and vector $v$, there is an eigenvalue $\lambda$ of $Y$ satisfying
% \[
% \vert\lambda - \theta\vert \le \Vert Y v - v \theta \Vert
% \]
%
% Notably, this result does not require that $(\theta, v)$ be a Ritz pair or have
% any variational structure whatsoever. Therefore, we can use
%
% \[
% \Delta \sigma_i = \sqrt{\Vert X X^T v_i - v_i \theta_i \Vert}
% \]
%
% as an error bar for using $\sqrt{\theta_i}$ as an estimate for a singular value,
% where $(\theta_i, v_i)$ is some eigenpair of $Q^T X X^T Q$ computed at any
% iteration of FlashPCA.

\section{Performance and accuracy of Lanczos vs. subspace iteration}

We compared FlashPCA and our code on a machine with 230 GB of RAM and two
Intel(R) Xeon(R) CPU E5-2697 v2 @ 2.70GHz processors, each with 12 cores.
We used the FlashPCA v1.2 release binary for Linux.
Our Julia code was run on a pre-release v0.5 version of Julia,
built with Intel Composer XE 2016.0.109 Linux edition,
and linked against Intel's Math Kernel Library v11.3.

Table~\ref{tab:runtime} shows a breakdown of the execution time in FlashPCA
compared to our Lanczos implementation on our synthetic genotype matrix.
We observe a large difference in the run time between our Lanczos implementation
and subspace iteration.
The bulk of the difference comes from Lanczos taking much fewer iterations
(as measured by matrix-vector products) to converge. For FlashPCA, most of the
run time is spent in an initial matrix multiply, $XX^T$, and the subsequent
subspace iterations.

\begin{table}
    \caption{Run time in seconds for a synthetic genotype matrix of
    size $41505\times81700$ generated using Algorithm~\ref{code:model}.
    FlashPCA 1: the published version of the FlashPCA program
    algorithm~\cite{abraham2014fast}, which normalizes but does not reorthogonalizes
    the vectors. FlashPCA 2: the current version of the FlashPCA program,
    which now by default normalizes and reorthogonalizes the vectors using dense
    QR factorization. This work: Golub-Kahan-Lanczos bidiagonalization without
    restarting, employing partial reorthogonalization.
    \label{tab:runtime}
    }
    \centering
    \begin{tabular}{lrrr}
                                   & FlashPCA 1 & FlashPCA 2   & this work     \\ \hline
        Preprocessing              &   54        &   61        & 54            \\
        Forming $XX^T$             &  971        &  926        & N/A           \\
        Subspace/Lanczos iteration &   41        & 1935        & 37            \\
        Postprocessing             &    7        &   11        & 0             \\ \hline
        Total                      & 1073        & 2933        & 81  \\ \hline
    \end{tabular}
\end{table}

Forming the matrix $XX^T$ is the default option in FlashPCA; while advantageous
for very wide matrices $\rho \approx 0$, it is an expensive step for our matrix,
which has $\rho \approx 2$.
Furthermore, we observe a performance penalty in using the Eigen linear algebra
library (used by FlashPCA) relative to MKL.
Whereas the former took 926 seconds to compute $XX^T$, the latter took only
306 seconds on the same machine, using the same number of threads.
This discrepancy may be system dependent.

\begin{figure}
    \centering
    \begin{subfigure}{0.72\textwidth}
        \begin{tikzpicture}
            \begin{semilogyaxis}[width=9.5cm, height=5cm, xmin = 0, xmax = 169, ymin=10e-15, ymax=10e1]
                \addplot[very thick] file {data/FlashPCAQRconv.txt};
            \end{semilogyaxis}
        \end{tikzpicture}
        \caption{FlashPCA 2}
    \end{subfigure}
    \begin{subfigure}{0.23\textwidth}
        \begin{tikzpicture}
            \begin{semilogyaxis}[width=3cm, height=5cm, ymin=10e-15, ymax=10e1]
                \addplot[very thick] file {data/JuliaPCAconv.txt};
            \end{semilogyaxis}
        \end{tikzpicture}
        \caption{this work}
    \end{subfigure}
    \caption{Convergence of the FlashPCA 2 and Lanczos-based algorithms. The
    convergence criterion used in FlashPCA 2 is \eqref{eq:flashpca-convergence},
    which in principle can be very different from the residual norm criterion used
    in our Lanczos-based implementation. Convergence of FlashPCA 2 is very slow,
    being nearly stagnant for many iterations before improving dramatically
    in the last iteration. In contrast, the Lanczos-based method shows essentially
    logarithmic convergence in the residual norm after the first few iterations.
    \label{fig:conv}
    }
\end{figure}
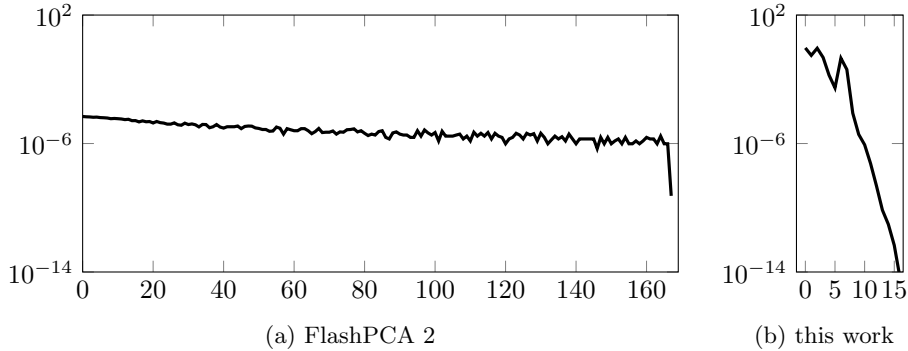

The subspace iterations form the bulk of the run time for FlashPCA 2. Even when
doing explicit reorthogonalization of the basis vectors, FlashPCA still performs
many more matrix-vector product-equivalents than our Lanczos implementation with
no restarting and partial reorthogonalization. Figure~\ref{fig:conv} contrasts
the convergence reported by FlashPCA 2 and our Lanczos implementation. Note that
the convergence criteria used by the two programs are very different---the former
uses the criterion \eqref{eq:flashpca-convergence}, which is very different from
the classical residual norm criteria commonly used by Lanczos-based
methods~\cite{Parlett1998}. Convergence of FlashPCA 2 is very slow,
being nearly stagnant for many iterations before improving dramatically
in the last iteration. In contrast, the Lanczos-based method shows essentially
logarithmic convergence in the residual norm after the first few iterations.

Finally, we compare in Table~\ref{tab:shootout} the relative performance
and accuracy of the subspace iteration methods against our Lanczos
implementations as well as two other established libraries,
ARPACK~\cite{Lehoucq1996} and PROPACK~\cite{Larsen1998}. ARPACK, while not
originally designed for singular value computations, can be used to
compute eigenvalues of the augmented matrix $\begin{pmatrix}0 & X\\X^T &
0\end{pmatrix}$, whose eigenvalues are the same as the singular values of
$X$. To measure performance, we count the number of matrix-vector product
(mvps) as well as the wall time using the threaded BLAS implementation of
MKL on 24 threads. To measure accuracy, we provide both
$\Vert\Delta\theta\Vert_1 = \sum_{j=1}^k \beta_n \vert U^{(n)}_{mj}
\vert$, the sum of all estimated errors for the singular values, as well
as the relative error of the 10th singular value
as computed by LAPACK. Our results show that our implementation of Lanczos
bidiagonalization provides qualitatively better singular values that the
subspace iteration methods provided by FlashPCA, even with QR
reorthogonalization. Furthermore, our performance is significantly better
than the standard tools, ARPACK and PROPACK.

\begin{table}
\begin{tabular}{|c|c|c|c|c|c|}
\hline
 & Algorithm & mvps & time (sec) & $\left\Vert \Delta\theta\right\Vert _{1}$ & rel. err. \tabularnewline
\hline
\hline
FlashPCA1 & Block power & 800 & 1073 & $2.55$ & $0.067$ \tabularnewline
\hline
FlashPCA2 & Block power & 33600 & 2933 & $0.0146$ & $1.88\cdot10^{-5}$ \tabularnewline
\hline
PROPACK & GKL (PRO, NR) & N/A & 120 & $5.2\cdot10^{-5}$ & $1.75\cdot10^{-11}$ \tabularnewline
\hline
ARPACK & GKL (FRO, IR) & 338 & 378 & $9.2\cdot10^{-6}$ & $8.34\cdot10^{-15}$ \tabularnewline
\hline
this work & GKL (FRO, TR) & 360 & 132 & $0.0037$ & $5.03\cdot10^{-8}$ \tabularnewline
\hline
this work & GKL (PRO, NR) & 190 & 81 & $2.5\cdot10^{-6}$ & $1.16\cdot10^{-14}$ \tabularnewline
\hline
\end{tabular}

\caption{Comparing various methods for computing the top 10 principal
components on a simulated genotype matrix of size $41505\times81700$.
Linear algebra kernels were run on MKL with 24 software threads except for
FlashPCA which ran on Eigen kernel.
Timings reported are wall times in seconds.
FlashPCA1 - FlashPCA with columnwise rescaling,
FlashPCA2 - FlashPCA with QR orthogonalization,
GKL - Golub--Kahan--Lanczos bidiagonalization,
PRO - partial reorthogonalization,
FRO - full reorthogonalization,
NR - no restart,
IR - implicit restart after 20 Lanczos vectors,
TR - thick restart after 20 Lanczos vectors,
mvps - number of matrix--vector products.
Termination criterion set to $\Vert \Delta Y\Vert = 10^{-8}$ for FlashPCA$n$,
and relative error in the 10th singular value compared with the value obtained from LAPACK.
\label{tab:shootout}
}
\end{table}

\section{Conclusions}

GWASs provide an exciting new data source for large scale matrix
computations, whose nominal dimensions are already on the order of $10^5
\times 10^5$ and will continue to grow rapidly in the near future.  The
statistical and computational demands of GWAS on genotype matrices
necessitate the best numerical algorithms and software.

We have implemented state of the art Lanczos bidiagonalization methods in
pure Julia, allowing us to compute the largest principal components of
genotype matrices more efficiently than any other tool currently being
used for genomics data analysis, and even outperforming some standard
packages for iterative eigenvalue and singular value computation such as
ARPACK and PROPACK. The implementation of these methods in the Julia
programming language provides a fast, practical software tool that permits
easy introspection into the inner workings of the Lanczos algorithms, as
well as experimentation into new methods with minimal fuss.

Further work may include generalizing the code to also handle block
Lanczos computations, which may further improve the performance of the
computation by making use of BLAS3 function calls. Imputation of missing
data will also become important in future data analysis, as nominal matrix
sizes grow and the number of incorrectly sequenced sites grows. We are
also implementing and studying into iterative methods for evaluating the
regression models used in GWASs.

\section*{Acknowledgments}

We thank the Julia development community for their contributions to free
and open source software. Julia is free software that can be downloaded
from \url{julialang.org/downloads}. The implementation of iterative SVD
described in this paper is available as the \verb|svdl| function in the
\href{https://github.com/JuliaLang/IterativeSolvers.jl}{IterativeSolvers.jl}
package. J.C. would also like to thank Jack Poulson (Stanford) and David
Silvester
(Manchester) for many insightful discussions.

\bibliography{bio,svd}

\end{document}